\documentclass[10pt,twoside]{art_large}
\usepackage{amssymb,amsbsy,amsmath,amsfonts,amssymb,amscd}
\usepackage[latin1]{inputenc}
\usepackage[english,francais]{babel}

\newtheorem{theorem}{Theorem}
\newtheorem{theorem-en}{Theorem}

\newtheorem{prop}[theorem]{Proposition}

\newtheorem{defn-en}[theorem-en]{Definition\rm}

\newcommand{\bbC}{\mathbb C}

 \ComParit{S0764-4442} \PIT{FLA} \PXMA{????} \Add{?}
 \Volume{???} \Year{2009} \FirstPage{1} \LastPage{4}
\AuteurCourant{Michaël Ayoul et Nguyen Tien Zung}
\TitreCourant{Galoisian obstructions to non-Hamiltonian integrability} \Journal
\Rubrique{Systèmes dynamiques}{Dynamical systems} \SousRubrique{}{}
\PresentePar{Prénom}{NOM}
\Recu{jour mois ann\'ee}{apr\`es r\'evision}{jour mois ann\'ee}
\begin{document}
\selectlanguage{english}
\title{Galoisian obstructions to non-Hamiltonian integrability}
\author{%
Michaël Ayoul and Nguyen Tien Zung}
\address{%
Institut de Mathématiques de Toulouse, UMR 5219 CNRS, Université Toulouse III \\
E-mail:  tienzung.nguyen@math.univ-toulouse.fr } \maketitle
\thispagestyle{empty}
%
%%%%%%%%%%%%%%%%%%%%%%%%%%%%%%%%%%%%%%%%%%%%%%%%%%%%%%%%%%%%
%%%  Abstract  %%%
%%%%%%%%%%%%%%%%%%
\selectlanguage{english}
\begin{Abstract}{%
We show that the main theorem of Morales--Ramis--Simo \cite{MRS-Galois2007} about Galoisian obstructions to meromorphic integrability of Hamiltonian systems can be naturally extended to the non-Hamiltonian case. Namely, if a dynamical system is meromorphically integrable in the non-Hamiltonian sense, then the differential Galois groups of the variational equations (of any order) along its solutions must be virtually Abelian. 
}\end{Abstract}
\begin{Resume}{%
Nous montrons la version non-hamiltonienne du théorème de
Morales-Ramis-Simo  \cite{MRS-Galois2007}. Plus précisément, si un système dynamique est méromorphiquement intégrable au sens non-hamiltonien, alors tous les groupes de Galois différentiels des équations variationels d'ordre arbitraire le long de ses solutions doivent \^etre virtuellement abéliens}\end{Resume}

\section{Introduction}

Recent work by Morales, Ramis and Simo \cite{MR-Galois2001,MRS-Galois2007} showed a very natural obstruction to meromorphic integrability (à la Liouville) of analytic Hamiltonian systems in terms of differential Galois theory.  In this note, we will show that their results can be naturally extended to the non-Hamiltonian case. The main idea is as follows: any dynamical system given by a vector field can be ``turned'' into a Hamiltonian system by a simple trick, namely the \emph{cotangent lifting}. If the original system is meromorphically integrable in the non-Hamiltonian sense, then the lifted system is meromorphically integrable à la Liouville, which implies that the differential Galois groups of the variational equations (of any order, along the non-stationary solutions) of the lifted system are virtually Abelian by Morales--Ramis--Simo theorem. On the other hand, the cotangent lifting process does not change the differential Galois groups of the vartiational equations at all. So we get the following result:

\begin{theorem} \label{thm:GaloisNonHam} 
Assume that a dynamical system given by a holomorphic vector field $X$ on a complex analytic variety $M$ is meromorphically integrable in the non-Hamiltonian sense, and let $\Gamma$ be an immersed complex curve (Riemann surface) in $M$ which is given by a non-stationary solution of $X$ (i.e. $X$ is tangent to $\Gamma$). Then for any narutal number $n \geq 1$, the differential Galois group of the variational equation of order $n$ of $X$ along $\Gamma$ is virtually Abelian (i.e., its Lie algebra is Abelian).
\end{theorem}

We hope that the above theorem will find many applications in non-Hamiltonian dynamics, e.g. for non-holonomic systems. We remark that some special cases of the above theorem were obtained earlier by Maciejewski, Przybylska, and Yoshida \cite{MP-Galois2008,MPY-Galois2008}. Our approach is different from theirs. 

In the rest of this note, we will explain the terminology used in the above theorem, and how this theorem can be deducted from the main theorem of  Morales--Ramis--Simo  \cite{MRS-Galois2007}.

\section{Non-Hamiltonian integrability and the cotangent lift}

For general dynamical (not necessarily Hamiltonian) systems, there is a natural notion of (non-Hamiltonian) integrability, which was probably studied first by Bogoyavlensky, but also independently by other people from different points of view, see e.g.  \cite{BaCu-Nonholonomic1999,Bogoyavlenskij-Integrability1998,Stolovitch1,Zung-Poincare2002,Zung-Torus2006}.  Let us recall this notion here:

Let $X$ be a vector field on a manifold $M$ of dimension $m$. Then we say that $X$ is integrable (in the non-Hamiltonian sense) if there are integers $k \geq 1, l \geq 0$ such that $k + l = m$, $p$ vector fields $X_1,\hdots,X_k$ on $M$ (which are linearly independent almost everywhere) which commute paire wise (i.e. $[X_i,X_j] = 0 \  \forall \ i,j$) and such that $X_1 = X$, and $q$ common first integrals $f_1,\hdots,f_l$ (which are functionally independent almost everywhere) for these vector fields on $M$ (i.e. $X_i(f_j) = 0  \  \forall \ i,j)$. If all the vector fields $X_i$ and common first integrals $f_j$ belong to some category, then we say that $X$ is integrable in that category. In particular, we will say that $X$ is \emph{meromorphically integrable} if the vector fields $X_i$ and the functions $f_j$ are meromorphic on $M$.

If a Hamiltonian system is integrable (à la Liouville, or in non-commutative sense) then it is also integrable in the non-Hamiltonian sense, though the inverse is not true, so non-Hamiltonian integrability is a weaker notion than Hamiltonian integrability (see, e.g.,  \cite{Zung-Poincare2002,Zung-Torus2006}).

Let us now explain the cotangent lift trick. If $\phi: M \to M$ is a diffeomorphism, then it gives rise to push-forward map (which is the inverse of the pull-back map) $\hat{\phi} := \phi_*  = (\phi^*)^{-1}: T^*M \to T^*M$ from the cotangent bundle of $M$ to itself, which projects to $\phi$ on $M$. We will say that $\hat{\phi}$ is the \emph{cotangent lift} of $\phi$. Similarly, if $X$ is a vector field on $M$, then it can naturally be lifted to a Hamiltonian vector field  on $T^*M$  (with respect to the canonical symplectic structure on  $T^*M$), which we will denote by $\hat{X}$. In fact, $X$ corresponds to a fiber-wise linear function $h_X$ on $T^*M$ via the pairing formula $h_X(p) = \langle p, X(x) \rangle  $ for $p \in T^*_xM$, and $\hat{X}$ is nothing but the Hamiltonian vector field of $h_X$.  We will say that  $\hat{X}$ is the \emph{cotangent lift of $X$}.  The vector field $\hat{X}$ contains $X$ as a sub vector field, via the identification of the zero section of the cotangent bundle $T^*M$ with $M$, and $\hat{X}$ projects to $X$ via the natural projection $\pi: T^*M \to M$. The flow of $\hat{X}$ on $T^*M$ is nothing but the cotangent lift of the flow of $X$ on $M$.

In terms of local coordinates, let $(z_1,\hdots,z_m)$ be a local coordinate system on $M$, and $(p_1,\hdots,p_m)$ be the dual coordinates on the cotangent fibers, i.e. if $p \in T^*_xM \subset T^*M$ has coordinates $z_1(p),\hdots,z_m(p),p_1(p),\hdots,p_m(p)$ then $z_i(p) = z_i(x)$ and $p$ is the covector $p = \sum_{i} p_i(p) d z_i $ over the point $x$.  The vector field $X$ has the form $X = \sum a_i {\partial \over \partial z_i}$, the corresponding fiber-wise linear function $h_X$ is $h_X = \sum a_i p_i$, or more precisely  $h_X (p) = \sum a_i (\pi(p)) p_i(p)$, and the cotangent lift $\hat{X}$ of $X$ is 
\begin{equation}
\hat{X} = \sum_{i} a_i {\partial \over \partial z_i} - \sum_i ({\sum_j {\partial a_j \over \partial z_i}} p_j) {\partial \over \partial p_i}.
\end{equation}.

If $\Phi (z) =(\Phi_1 (z),\hdots,\Phi_m (z)) : \bbC^m \to \bbC^m$ is a local diffeomorphism, $z = (z_1,\hdots,z_m)$, then its cotangent lift $\hat{\Phi} : T^*\bbC^m \to T^*\bbC^m$, then $\hat{\Phi}$ in the coordinate system $(z,p) = (z_1,\hdots,z_m,p_1,\hdots,p_m)$ has the following form: 
\begin{equation}\label{eqn:1}
\hat{\Phi} (z,p) = (\Phi_1(z),\hdots,\Phi_m(z),\Psi_1(z,p),\hdots,\Psi_m(z,p))
\end{equation}
where
\begin{equation} \label{eqn:2}
(\Psi_1(z,p),\hdots,\Psi_m(z,p))^T = \left( {\partial \Phi_j \over \partial z_i} \right)_{i,j=1,\hdots,m}^{-1} (p_1,\hdots,p_m)^T .
\end{equation}

The cotangent lift relates many notions in Hamiltonian dynamics with notions in non-Hamiltonian dynamics. For example, one can verify easily that if $X$ is a vector field which vanishes at the origin and which is in Poincaré-Dulac normal form, then its cotangent lift is in Birkhoff normal form in the corresponding coordinates. In this note, we will need the following very simple but useful result:

\begin{prop} \label{prop:CotangentLiftIntegrable} If $X$ is (meromorphically) integrable in the non-Hamiltonian sense, then its cotagent lift $\hat{X}$ is (meromorphically) integrable à la Liouville.
\end{prop}

\emph{Proof.} Let $X_1,\hdots,X_k$ ($X_1 = X$) and $f_1,\hdots,f_l$ ($k+l = m$) be commuting vector fields and common first integrals which provide the non-Hamiltonian integrability of $X$. Consider the following $m$ functions $H_1,\hdots,H_m$ on the symplectic manifold $T^*M$: $H_1 = h_{X_1},\hdots,H_k = h_{X_k},H_{k+1} = \pi^* f_1,\hdots, H_m = \pi^* f_l$. Then $\hat{X}$ is the Hamiltonian vector field of $H_1$, and the functions  $H_1,\hdots,H_m$ Poisson-commute pairwise. Indeed, $[X_i,X_j] = 0$ implies that $\{ H_i,H_j\} = 0$ (for $i,j \leq k$); $X_i(f_j) = 0$ implies that $\{H_i, H_{k+j}\} = 0$, and we always have $\{\pi^*f,\pi^*g\} = 0$ for any two functions $f,g$ on $M$. Thus $\hat{X} = X_{H_1}$ is an integrable Hamiltonian system à la Liouville with the momentum map $(H_1,\hdots,H_m)$.  {\hfill $\Box$}

\section{Variational equations and differential Galois groups}

Let us birefly recall how to construct the variational equation of order $n$ ($n \in {\mathbb N}$) along a non-stationary solution  of a dynamical system given by a complex analytic vector field $X$, and the corresponding differential Galois group. For details see \cite{MRS-Galois2007}. We will assume that $\Gamma$ is an immersed Riemann surface in the manifold $M$, which is tangent to $X$ ($X$ is non-trivial on $\Gamma$ but may vanish et some points on $\Gamma$). Denote by ${\mathcal M}(\Gamma)$ the field of meromorphic functions on $\Gamma$. The vector field $X$ (resitricted to $\Gamma$) turns ${\mathcal M} (\Gamma)$ into a differential field, with the differential $\partial$ given by $X$. For each point $x \in \Gamma$ denote by $J^n_x$ the set of $n$-jets of holomorphic functions in a neighborhood of $x$ in $M$. Then $J^n_x$ is a finite-dimensional vector space, and the union $J^n_\Gamma = \bigcup_{x \in \Gamma} J^n_x$ is a locally trivial complex vector bundle over $\Gamma$. The vector field $X$ induces a natural linear connection on $J^n_\Gamma$ in the following way: for each pair of points $x, x' \in \Gamma$ and a path $\alpha$ in $\Gamma$ connecting $x$ to $x'$, there is a unique path $\alpha_\bbC$ in $\bbC$ from $0$ to some complex number such that the integration of $X$ along this path $\alpha_\bbC$ of complex times will move $x$ to $x'$. Let $\phi_{\alpha}$  denote the germ (near $x$) of the integration of $X$ along the path $\alpha_\bbC$ of complex times. Then $\phi_{\alpha}$ is a germ of complex diffeomorphism from a neighborhood of $x$ in $M$ to a neighborhood of $x'$. The $n$-jet of the push-forward of  $\phi_{\alpha}$ is a linear isomorphism from $J^n_x$ to $J^n_{x'}$, which we will denote by $J^n_\alpha(X) : J^n_x \to J^n_{x'}$. This linear isomorphism depends only on the homotopy class of $\alpha$, and is the linear transport map of our linear connection on $J^n_\Gamma$. This linear connection on $J^n_\Gamma$ is what is called the \emph{dual variational equation of order $n$ along $\Gamma$} in \cite{MRS-Galois2007}. In order to view the above linear connection on $J^n_\Gamma$ as a system of linear differential equations on $\Gamma$, one needs to trivialize $J^n_\Gamma$ meromorphically, i.e. fix a basis of meromorphic sections on $J^n_\Gamma$ (such a meromorphic trivialization exists and is unique up to meromorphic gauge transformations, see e.g. the appendices in \cite{MR-Galois2001,MRS-Galois2007}.) Denote the Picard--Vessiot extension of ${\mathcal M}(\Gamma)$ associated to the corresponding linear system of equations on  $J^n_\Gamma$ by ${\mathcal M}^n_X(\Gamma)$ (it is generated by the coefficients of the fundamental solution of the system of linear differential equations, and does not depend on the choice of the meromorphic trivialization.)  We will denote by $Gal^n_X( \Gamma)$ the differential Galois group $Gal( {\mathcal M}^n_X(\Gamma)/ {\mathcal M}(\Gamma))$ of the extension ${\mathcal M}^n_X(\Gamma)/ {\mathcal M}(\Gamma)$, i.e. the group of differential automorphisms of ${\mathcal M}^n_X(\Gamma)$ which are identity on  ${\mathcal M}(\Gamma)$. Recall that these groups are linear algebraic groups. With the above notations, the main result  of Morales--Ramis--Simo \cite{MRS-Galois2007} may be formulated as follows:

\begin{theorem} \label{thm:MRS}
  If $M$ is a symplectic manifold and $X$ is a Hamiltonian vector field which is meromorphically integrable à la Liouville (i.e. there is a set of functionally independent meromorphic first integrals in involution, whose cardinal is half the dimension of $M$), then the groups  ${\mathcal G}^n_X(\Gamma)$ are virtually Abelian.
\end{theorem}

\section{Cotangent lift and differential Galois groups}

We keep the notations of the previous sections. We will view $X$ as a sub vector field of its cotangent lift $\hat{X}$, via the identification of $M$ with the zero section of $T^*M$. Then $\Gamma$ is also an integral curve of $\hat{X}$.

\begin{theorem}  \label{thm:CotangentLiftGalois}
We have $Gal^n_{\hat{X}}(\Gamma) \cong Gal^n_{X}(\Gamma)$ (for any holomorphic vector field $X$ on a complex manifold $M$, integral curve $\Gamma$ of $X$ in $M$, and natural number $n$).
\end{theorem}

As we already mentioned in the introduction, the proof of Theorem \ref{thm:GaloisNonHam}  is a direct consequence of Proposition \ref{prop:CotangentLiftIntegrable}, Theorem \ref{thm:MRS}, and Theorem \ref{thm:CotangentLiftGalois}. So in order to prove Theorem  \ref{thm:GaloisNonHam}, it remains to prove Theorem \ref{thm:CotangentLiftGalois}.

\emph{Proof of Theorem \ref{thm:CotangentLiftGalois}.}  The main idea is that, by making the cotangent lift, we don't create any ``essentially new'' information (all the information is already contained in the original system). In particular, we don't create any ``essentially new'' function, so the Picard--Vessiot extensions don't change.

For simplicity, let us assume for the moment that for almost every $x \in \Gamma$ there is a local holomorphic coordinate system $(z_1^x, \hdots, z_m^x)$ in $M$ centered at $x$, such that the family $(z_1^x, \hdots, z_m^x)$ depends on $x$ in a meromorphic way. (In other words, the functions $z_i^x$ may be viewed as meromorphic functions defined near the diagonal in $\Gamma \times M$, and vanish on the diagonal). This meromorphic family of local coordinate systems will give a meromorphic trivialization of $J^n_\Gamma$ via the coefficients of Taylor expansions. Let us fix a regular point $x_0 \in \Gamma$. For $x \in \Gamma$ near $x_0$, the ``time $x-x_0$'' flow of $X$ moves a neighborhood of $x_0$   to a neighborhood of $x$  by a local diffeomorphism $\phi_X^{x_0,x}$. With respect to the local coordinate systems  $(z_1^{x_0}, \hdots, z_m^{x_0})$ and  $(z_1^x, \hdots, z_m^x)$, $\phi_X^{x_0,x}$ becomes a local diffeomorphism from $(\bbC^m,0)$ to $(\bbC^m,0)$, which we will denote by $\Phi_x^X = (\Phi_{1,x}^X,\hdots, \Phi_{m,x}^X): (\bbC^m,0) \to (\bbC^m,0)$, where each $\Phi_{i,x}^X = \Phi_{i,x}^X(z_1,\hdots,z_m)$ is a local function of $m$ variables $(z_1,\hdots,z_m)$ which depend on the parameter $x$.  The value at 0 of the partial derivatives of the functions $\Phi_{i,x}$ (with respect to the variables $(z_1,\hdots,z_m)$) up to order $n$  give us (germs of) local holomorphic functions in $x$ in a neighborhood of $x_0$, and the Picard--Vessiot extension ${\mathcal M}^n_X(\Gamma)$ of the variational equation of order $n$ is nothing but the differential extension of $\mathcal{M}(\Gamma)$ by these functions (i.e. values of the partial derivatives up to order $n$ of the components of the diffeomorphism $\Phi_x$ at $z=0$). 

Each local coordinate system $(z_1^x, \hdots, z_m^x)$ gives rise to a coresponding set of dual coordinates $(p_1^x,\hdots,p_m^x)$ in the cotangent bundle $T^*M$ (near  $T^*_xM$). Thus we have a meromorphic family $(z_1^x, \hdots, z_m^x, p_1^x,\hdots,p_m^x)$  of local coordinate systems in $T^*M$ near the points of $\Gamma$ ($\Gamma \subset M$, $M$ being identified with the zero section of $T^*M$.) With respect to these local coordinate systems,  the flow map $\Phi_x^{\hat{X}}$ of the lifted vector field $\hat{X}$ is the cotangent lift of the original flow map $\Phi_x^{X}$.  It follows from Formulas (\ref{eqn:1}) and (\ref{eqn:2}) that the values at $z=p=0$ of the partial derivatives up to (total) order $n$ in $z$ and $p$ of the components of $\Phi_x^{\hat{X}}$ don't give anything new (i.e. anything which does not belong to ${\mathcal M}^n_X(\Gamma)$) compared to the values at $z=0$ of the partial derivatives up to order $n$ of the components of the diffeomorphism $\Phi_x$. It implies that ${\mathcal M}^n_{\hat{X}}(\Gamma)$ is equal to ${\mathcal M}^n_X(\Gamma)$, which implies that $Gal^n_{\hat{X}}(\Gamma) \cong Gal^n_{X}(\Gamma)$.

In general, we don't know for sure if we will be able to find a meromorphic family of loal coordinate systems $(z_1^x, \hdots, z_m^x)$. However, it follows from the meromorphic trivialisability of $J^n_\Gamma$ that we can choose $(z_1^x, \hdots, z_m^x)$ so that their $n$-jets depend meromorphically on $x$ (though maybe the functions $z_i^x$ themself are not meromorphic in $x$). But since we are in fact working with $n$-jets only, the above arguments are still valid, and we still have $Gal^n_{\hat{X}}(\Gamma) \cong Gal^n_{X}(\Gamma)$.  {\hfill $\Box$}

\vspace{0.5cm}

\emph{Acknowledgements}. This work is partially supported by the ANR project ``GIMP'' (ANR-05-BLAN-0029-01) and the ANR project ``Integrabilité réelle et complexe en mécanique hamiltonienne'' (JC05-41465). The authors would like to thank Jean Pierre Ramis for explaining to them several aspects of the papers \cite{MR-Galois2001,MRS-Galois2007}. They would also like to thank the members and the guests of the ANR project ``'Intégrabilité ...'', especially Jacques-Arthur Weil and Guy Casale, for their interest in this work and helpful remarks.

%\bibliographystyle{amsplain}
%\bibliography{zung}
%\end{document}
\providecommand{\bysame}{\leavevmode\hbox to3em{\hrulefill}\thinspace}

\end{document}